\documentclass[MSNbibl,number,citesort,dvips]{arxbj}
\usepackage{upgreek}
\usepackage{dcolumn}
\usepackage{graphicx}

\aid{0}
\volume{19}
\issue{5A}
\pubyear{2013}
\firstpage{1839}
\lastpage{1854}
\doi{10.3150/12-BEJ432} 

\makeatletter

\newcolumntype{d}[1]{D{.}{.}{#1}}
\def\eqref#1{(\ref{#1})}
\newtheorem{lemma}{Lemma}
\newtheorem{theorem}{Theorem}
\newcommand{\bfh}{\mathbf{h}}
\newcommand{\bfy}{\mathbf{y}}
\newcommand{\bfd}{\mathbf{d}}
\newcommand{\bfg}{\mathbf{g}}
\newcommand{\bfs}{\mathbf{s}}
\newcommand{\bfmu} {\bolds{\mu}}
\newcommand{\bfbeta}{\bolds{\beta}}
\newcommand{\bfvarepsilon} {\bolds{\varepsilon}}
\newcommand{\bfepsilon}{\bolds{\epsilon}}
\def\be{\begin{eqnarray}}
\def\ee{\end{eqnarray}}
\def\bse{\begin{eqnarray*}}
\def\ese{\end{eqnarray*}}

\makeatother

\begin{document}
\begin{frontmatter}

\title{Optimal variance estimation without estimating the mean function}
\runtitle{Optimal variance estimation}

\begin{aug}
\author[1]{\fnms{Tiejun} \snm{Tong}\corref{}\thanksref{1}\ead[label=e1]{tongt@hkbu.edu.hk}},
\author[2]{\fnms{Yanyuan} \snm{Ma}\thanksref{2}\ead[label=e2,text=ma@stat. tamu.edu]{ma@stat.tamu.edu}} \and
\author[3]{\fnms{Yuedong} \snm{Wang}\thanksref{3}\ead[label=e3]{yuedong@pstat.ucsb.edu}}
\runauthor{T. Tong, Y. Ma and Y. Wang} 
\address[1]{Department of Mathematics, Hong Kong Baptist University,
Hong Kong.\\ \printead{e1}}
\address[2]{Department of Statistics, Texas A\&M University, College
Station, TX 77843, USA.\\ \printead{e2}}
\address[3]{Department of Statistics and Applied Probability,
University of California, Santa Barbara, CA 93106, USA. \mbox{\printead{e3}}}
\end{aug}

\received{\smonth{8} \syear{2011}}
\revised{\smonth{12} \syear{2011}}

%
\begin{abstract}
We study the least squares estimator in the residual variance
estimation context. We show that the mean squared
differences of paired observations are asymptotically normally
distributed. We
further establish that, by regressing the mean squared differences of
these paired observations on the squared distances between paired
covariates via a simple least squares procedure, the resulting variance
estimator is not only asymptotically normal and root-$n$ consistent,
but also reaches the optimal bound in terms of estimation variance.
We also demonstrate the advantage of the least squares estimator in
comparison with existing methods in terms of the
second order asymptotic properties.
\end{abstract}

%

\end{frontmatter}
%

\section{Introduction}\label{sec1}

Consider the following nonparametric regression model
\begin{eqnarray}\label{model}
y_{i}=g(x_i)+\varepsilon_{i},\qquad 0 \le x_i \le1, i=1,\dots,n,
\end{eqnarray}
where $y_i$ is the observation of the mean function $g$ evaluated
at design point $x_i$ plus random error $\varepsilon_{i}$. We assume
that $\varepsilon_{i}$'s are independent and identically distributed
with mean zero and variance $\sigma^2$. Many nonparametric regression
methods have been developed to estimate the mean function $g$ in the
literature. Often, for choosing the amount of smoothing, testing
goodness of fit or estimating model complexity, one needs an estimate
of $\sigma^2$ that does not require estimating the mean function $g$
first \cite{Eubank1990,Gasser1991,Ye1998}. For example, an estimate of
$\sigma^2$ is required in the unbiased risk criterion for selecting
the smoothing parameter in spline smoothing (see Section 3.3 in
\cite{Wangbk}).

One popular class of estimators of $\sigma^2$ which bypasses the
estimation of $g$ is the so-called
difference-based estimators. The basic idea of difference-based
estimation is to use differences to remove trend in the mean
function. Assume that $0 \leq x_1 \leq\cdots\leq x_n \leq1$.
Rice \cite{Rice1984} proposed the first-order difference-based
estimator
\begin{equation}
{\hat\sigma_R^2}= \frac{1}{2(n-1)} \sum_{i=2}^{n} (y_i-y_{i-1})^2 .
\end{equation}
Gasser, Sroka and Jennen-Steinmetz \cite{Gasser1986} and Hall, Kay and
Titterington \cite{HallKayTit1990} extended the
Rice estimator to the second- and higher-order difference-based
estimators, respectively. More difference-based estimators can be
found in \cite{Dette1998,Muller2003}.

Tong and Wang \cite{Tong2005} proposed a variation of the
difference-based estimator.
For simplicity, consider the equally-spaced design where $x_i=i/n$.
Define the lag-$k$ Rice estimators as
\begin{eqnarray}\label{lag.Rice}
s_k = \frac1 {2(n-k)} \sum_{i=1}^{n-k} {(y_{i+k}-y_i)^2},
\qquad k=1,2,\dots.
\end{eqnarray}
For any $k=\mathrm{o}(n)$, it can be shown that $\mathrm{E}(s_k)= \sigma
^2 + J
d_k + \mathrm{o}(d_k)$ where $J= \int_0^1 \{ g'(x) \}^2\,\mathrm{d}x/2$ and
$d_k=k^2/n^2$. That is, the lag-$k$ Rice estimator overestimates
$\sigma^2$ by $J d_k$. To reduce bias, they proposed
fitting a linear regression model
\begin{eqnarray}\label{fitted.model}
s_k = \beta_0 + \beta_1 d_k + \epsilon_k,\qquad k=1,\dots,m,
\end{eqnarray}
where $m=\mathrm{o}(n)$ and using the least squares type of estimate of the
intercept as
an estimate of~$\sigma^2$.

For ease of notation, let $\bfs=(s_1,\dots,s_m)^T$, $\bfbeta=(\beta
_0,\beta_1)^T$, $\bfepsilon=(\epsilon_1,\dots,\epsilon_m)^T$,
$\mathbf{1}=(1,\dots,1)^T$, $\bfd=(d_1,\dots,d_m)^T$ and $X=(\mathbf{1},
\bfd)$ be the design matrix.
Then (\ref{fitted.model}) leads to $\bfs= X\bfbeta+ \bfepsilon$.
Note that $s_k$ is the average of $(n-k)$ lag-$k$ differences and there
are a total of $N=(n-1) + (n-2) + \cdots+ (n-m) = nm - m(m+1)/2$ pairs
of differences involved in the regression.
Tong and Wang \cite{Tong2005} assigned weight $w_k=(n-k)/N$ to the observation
$s_k$ and then fitted the linear regression using the weighted least
squares with weight matrix $W=\operatorname{diag}(w_1,\dots,w_m)$.
This results in $\hat{\bfbeta}_{\mathrm{WLS}} = (X^T W^{-1}X)^{-1}
X^T W^{-1}\bfs$. Consequently,
the weighted least squares estimator of $\sigma^2$ is
\begin{eqnarray}\label{est1}
\hat\sigma^2 = \hat\beta_{0,{\mathrm{WLS}}} = \sum_{k=1}^m
w_ks_k -
\hat\beta_{1,{\mathrm{WLS}}} \bar d_w,
\end{eqnarray}
where $\bar d_w=\sum_{k=1}^m w_kd_k$ and $\hat
\beta_{1,{\mathrm{WLS}}} = {{\sum_{k=1}^m w_ks_k(d_k-\bar
d_w)}/{\sum_{k=1}^m w_k(d_k-\bar d_w)^2}}$.
For simplicity, the above weighted least squares estimator
$\hat\sigma^2$ is referred to as the least squares estimator in this
paper. In Section \ref{sec3}, we will show that the above weighted
least squares
estimator is asymptotically equivalent to the ordinary least squares
estimator and the generalized least squares estimator which takes into
account the correlations between~$s_k$'s.

In this paper, we investigate the asymptotic distribution and
efficiency of the least squares estimator. We show that the least
squares estimator is asymptotically normally distributed in Section \ref{sec2}.
We further show that the least squares estimator is asymptotically
equivalent to the generalized least squares
estimator where correlations among $s_k$ are accounted for in Section
\ref{sec3}. In Section \ref{sec4}, we derive the optimal efficiency
bound for any
estimation procedure and show that
the least squares estimator reaches this optimal efficiency bound.
In Section \ref{sec5}, we derived the mean squared error (MSE) for M{\"u}ller and Stadtm{\"u}ller's \cite{Muller1999} estimator
and then compare it to the least squares estimator.
A~real example is also provided.
Finally, we conclude the paper in Section \ref{sec6} with some
simulation studies.

\section{Least squares estimator}\label{sec2}

Let $\bfy= (y_1,\dots,y_n)^T$, $\bfg=(g(x_1),\dots,g(x_n))^T$ and
$\bfvarepsilon= (\varepsilon_1,\dots,\varepsilon_n)^T$.
Then $\bfy=\bfg+\bfvarepsilon$. Let
$\gamma_i=\mathrm{E}(\varepsilon^i)/\sigma^i$ for $i=3,4$, and
$\stackrel{\mathcal{D}}\to$ denote convergence in distribution.
Assume that $\gamma_4 > 1$.
We first establish asymptotic normality for the Rice estimator.

\begin{theorem} \label{theorem1}
Assume that $g$ has a bounded second derivative. For any $k=n^r$ with
$0<r<3/4$, the lag-$k$ Rice estimator satisfies
$\sqrt{n}(s_k -\sigma^2) \stackrel{\mathcal{D}}\to N(0,
\gamma_4\sigma^4)$ as $n\to\infty$.
\end{theorem}

Proof of Theorem \ref{theorem1} can be found in Appendix \ref{appA}.
Next, we establish
asymptotic normality for the least squares estimator (\ref{est1}).
Following the result in \cite{Tong2005}, the least squares
estimator (\ref{est1}) has a quadratic form $\hat\sigma^2 =
\bfy^TD\bfy/\operatorname{tr}(D)$,
where $D=(d_{ij})_{n\times n}$ is a symmetric matrix with elements
\[
d_{ij}=\cases{\displaystyle
\sum_{k=1}^m b_k + \sum_{k=0}^{\min(i-1,n-i,m)} b_k, &\quad $1 \le i=j \le n$,
\vspace*{2pt}\cr
-b_{|i-j|}, &\quad $0<|i-j|\le m$, \vspace*{2pt}\cr
0, & \quad otherwise,
}
\]
where $b_0=b_{m+1}=0$ and
$b_k=1-{\bar d_w (d_k-\bar d_w)/\sum_{k=1}^m w_k(d_k-\bar d_w)^2}$
for $k=1,\dots,m$.

\begin{theorem} \label{theorem2}
Assume that $g$ has a bounded second derivative and $\mathrm
{E}(\varepsilon^6)$ is finite.
Then for any $m=n^r$ with $0<r<1/2$, the least squares estimator
$\hat\sigma^2$ satisfies
$\sqrt{n}(\hat\sigma^2 - \sigma^2) \stackrel{\mathcal{D}}\to
N\{0, (\gamma_4-1)\sigma^4\}$ as $n\to\infty$.
\end{theorem}

Proof of Theorem \ref{theorem2} can be found in Appendix \ref{appB}.
Given that $\mathrm{E}(\varepsilon^6)$ is finite, Theorems \ref
{theorem1} and~\ref{theorem2} show that the least squares estimator is
more efficient than the Rice estimator.
Theorem \ref{theorem2} also indicates that the least squares estimator
is as efficient as the
sample variance based on independent and identically distributed
samples, regardless of whether the
unknown mean function is a constant or not.

Theorem \ref{theorem2} can be used to construct confidence intervals
for $\sigma^2$. Assume that $n > (\gamma_4-1)
z_{\alpha/2}^2$ where $z_\alpha$ is the upper $\alpha$th percentile
of the standard normal distribution. Then an approximate $1-\alpha$
confidence interval for $\sigma^2$ is $[\hat\sigma^2/\{1+z_{\alpha
/2}\sqrt{(\gamma_4-1)/n}\}, \hat\sigma^2/\{1-z_{\alpha/2}\sqrt
{(\gamma_4-1)/n}\}]$.
For the special case when the $\varepsilon_i$'s are distributed from
$N(0,\sigma^2)$, we have $\gamma_4=3$.
In general, the parameter $\gamma_4$ can be replaced by an estimate.
Finally, by Box \cite{Box1954} and Rotar \cite{Rotar1973},
the finite sample distribution of $\hat\sigma^2$ can be approximated
by the scaled chi-squared distribution, $(\sigma^2/\nu)\chi^2(\nu)$,
where $\nu= \{ \operatorname{tr}(D) \}^2/\operatorname{tr}(D^2)$.

\section{Generalized least squares estimator}\label{sec3}

In Appendix \ref{appC}, we show that, for any $1\leq b<k=n^r$ with $0<r<2/3$,
$\operatorname{Cov}(s_b,s_k) = n^{-1}(\gamma_4-1)\sigma^4 + \mathrm{o}(n^{-1})$.
Combined with the results in Theorems \ref{theorem1}, we have
$\operatorname{Corr}(s_b,s_k)
\to
(\gamma_4-1)/\gamma_4$ as $n\to\infty$. In the case when the
$\varepsilon_i$'s are normally distributed, $\gamma_4=3$ and the
correlation coefficients between the lag-$k$ Rice estimators are all
asymptotically equal to~$2/3$.

In the construction of the least squares estimator in Section
\ref{sec2}, we have ignored the correlation between $s_k$'s. Given that the
correlation between lag-$k$ Rice estimators are high, a
natural question is whether the least squares estimator can be
improved by the following generalized least squares estimator
\begin{eqnarray}\label{GLS.est}
\hat{\bfbeta}_{\mathrm{GLS}} =
(X^T\Sigma^{-1}X)^{-1}X^T\Sigma^{-1}\bfs,
\end{eqnarray}
where $\Sigma=\gamma_4\sigma^4 \{ (1-\rho) I + \rho\mathbf{1}^T \mathbf{
1} \}/n$ is the asymptotic variance--covariance matrix,
$\rho=(\gamma_4-1)/\gamma_4$, and $I$ is the identity matrix.
It is known that $\hat{\bfbeta}_{\mathrm{GLS}}$ is the best linear
unbiased estimator of $\bfbeta$ \cite{Kariya2004bk}. Since $\Sigma$
has the compound symmetry structure and the first column of $X$ is
$\mathbf{1}$, by McElroy \cite{McElroy1967}, the generalized least squares
estimator $\hat{\bfbeta}_{\mathrm{GLS}}$ is identical to the
ordinary least squares estimator
$\hat{\bfbeta}_{\mathrm{OLS}} = (X^TX)^{-1} X^T\bfs$.
Furthermore, for any $m=\mathrm{o}(n)$, it is not difficult to show that
$\hat{\bfbeta}_{\mathrm{WLS}}$
is equivalent to $\hat{\bfbeta}_{\mathrm{OLS}}$.
Therefore, $\hat{\bfbeta}_{\mathrm{OLS}}$, $\hat{\bfbeta}_{\mathrm
{GLS}}$ and
$\hat{\bfbeta}_{\mathrm{WLS}}$ are all asymptotically equivalent.

\section{\texorpdfstring{The optimal efficiency bound for estimating $\sigma^2$}
{The optimal efficiency bound for estimating sigma2}}\label{sec4}

In this section, we derive the optimal semiparametric efficiency bound
for estimating $\sigma^2$ in model~\eqref{model} for any estimation
procedure and show that the least squares estimator reaches this bound.

Consider the estimation of $\sigma^2$ in model (\ref{model})
regardless of how the estimation is carried out. For simplicity, we omit
the subindex $i$. Under (\ref{model}), the only assumption is that
$\varepsilon=Y-g(X)$ are independent and identically distributed with
mean zero, and are independent of $X$. Denote the model of the
probability density function of $\varepsilon$ as $\eta(\varepsilon)$.

The probability density function model of $(x, y)$ can be written as
$f_X(x)\eta\{y-g(x)\}=f_X(x)\eta(\varepsilon)$, where $f_X(\cdot)$
is a
marginal probability density function model of $X$ and $\eta$ is a
probability density function model that ensures
zero mean, i.e.,
$\int\eta(\varepsilon)\,\mathrm{d}\varepsilon=1$ and $\int\varepsilon
\eta
(\varepsilon)\, \mathrm{d}\varepsilon=0$.
Viewing $f_X,
\eta$ and $g$ as the nuisance parameters and $\sigma^2=\mathrm
{E}(\varepsilon^2)$ as the parameter of
interest, this becomes a semiparametric problem and one can
derive the efficient influence function through projecting any
influence function
onto the tangent space associated with $f_X$, $\eta$ and $g$.

Simple calculation
yields the tangent space of model (\ref{model}) to be
\be\label{eq:lambda}
\Lambda_{\mathcal{T}}&=&
\{ h(x)+f(\varepsilon)+\eta_0'(\varepsilon)/\eta
_0(\varepsilon)a(x)\dvtx\nonumber\\[-8pt]\\[-8pt]
&&\hspace*{4pt} \forall h, f \mbox{ such
that } \mathrm{E}(h)=0,
\mathrm{E}(f)=\mathrm{E}(\varepsilon f)=0, \mbox{and } \forall a
\},\nonumber
\ee
where $\eta_0(\cdot)$ denotes the true probability density function
of $\varepsilon$.
Following the procedure in Chapter 4 of \cite{Tsiatis2006}, we
consider an arbitrary parametric submodel, denoted as $\eta
(\varepsilon, \bfmu)$. Here
$\bfmu$ is a finite dimensional vector of parameters and
there exists $\bfmu_0$,
such that $\eta(\varepsilon,\bfmu_0)=\eta_0(\varepsilon)$. In addition,
$\eta(\varepsilon,\bfmu)$ is a valid probability density function and
$\int\varepsilon\eta(\varepsilon;\bfmu)\,\mathrm{d}\varepsilon=0$ for all
$\bfmu$ in a local
neighborhood of $\bfmu_0$. We have
$
\partial\int\varepsilon^2 \eta(\varepsilon, \bfmu)\,\mathrm
{d}\varepsilon
/\partial\bfmu
=\mathrm{E}(\varepsilon^2S_{\bfmu})$,
where $S_{\bfmu}=\partial\log\eta(\varepsilon, \bfmu)/\partial
\bfmu$
is the score vector with respect to $\bfmu$. Hence,
$\varepsilon^2-\sigma^2$ is a valid influence function.
We decompose $\varepsilon^2-\sigma^2$ into
\[
\varepsilon^2-\sigma^2
=\{\varepsilon^2-\sigma^2+\gamma_3\sigma^3\eta_0'(\varepsilon
)/\eta_0(\varepsilon)\}
-\gamma_3\sigma^3\eta_0'(\varepsilon)/\eta_0(\varepsilon).
\]
It is not difficult to verify that
$\varepsilon^2-\sigma^2+\gamma_3\sigma^3\eta_0'(\varepsilon)/\eta
_0(\varepsilon)$
satisfies the requirement on $f$ in (\ref{eq:lambda}). Hence, it is a
qualified $f(\varepsilon)$ function. Letting
$a(x)$ in (\ref{eq:lambda}) be $-\gamma_3\sigma^3$ yields
$-\gamma_3\sigma^3\eta_0'(\varepsilon)/\eta_0(\varepsilon)$.
Thus, $\varepsilon^2-\sigma^2\in\Lambda_{\mathcal{T}}$, and consequently
it is the efficient influence function.
The corresponding efficient estimation variance is
$n^{-1}\mathrm{E}\{(\varepsilon^2-\sigma^2)^2\}=n^{-1}(\gamma_4-1)\sigma
^4$, which
agrees with the result in Theorem \ref{theorem2}. This shows that the
least squares estimator is indeed optimal in terms of its
estimation variability among the class of all root-$n$ consistent
estimators.

In the above derivation, we have not taken into account
that $X_i$'s are actually equally spaced instead of being
random. However, assuming $f_X(x)$ to be uniform or more generally
assuming $f_X(x)$ to have any particular form does not change the
efficiency result. This is because the calculation
relies on the property of $\varepsilon$ only, which is
independent of~$X$.

\section{Variance estimator of M{\"{u}}ller and Stadtm{\"{u}}ller}\label{sec5}

M{\"u}ller and Stadtm{\"u}ller \cite{Muller1999} proposed a similar
least squares type
estimator for the equally-spaced design where $x_i=i/n$.
Define
\begin{eqnarray*}
z_k = \frac{1}{2(n-L)}\sum_{i=1}^{n-L}(y_{i+k}-y_i)^2,\qquad 1\leq
k\leq L,
\end{eqnarray*}
where $L=L(n)\geq1$. In the context of testing if the mean function
contains jump discontinuities, M{\"u}ller and Stadtm{\"u}ller \cite
{Muller1999} fitted a
linear model that regresses $z_k$ on two independent variables,
one for the sum of the squared jump sizes and the other for the
integrated squared first derivative, and then estimate the residual
variance as the intercept. In the case when the function is smooth,
that is, when the sum of the squared jump sizes equals to zero, the
variance estimator in \cite{Muller1999} reduces to
\begin{eqnarray}\label{est.muller}
\hat\sigma^2_{\mathrm{MS}}
= {3\over L(L-1)(L-2)}\sum_{k=1}^L \{ 3L^2+3L+2-6(2L+1)k+10k^2 \} z_k.
\end{eqnarray}

The dependent variable $z_k$ in \cite{Muller1999} uses the first
$n-L$ terms in the lag-$k$ Rice estimator $s_k$ while the last $L-k$
terms are ignored. This makes $z_k$ a less efficient estimator of
$\sigma^2$, especially when $L-k$ is large. In addition, noting that
$\hat\sigma^2_{\mathrm{MS}}$ is a weighted average of $z_k$ with
larger weights assigned to small $k$ and more terms are ignored with
small $k$, the efficiency loss of $\hat\sigma^2_{\mathrm{MS}}$ over
$\hat\sigma^2$ can be severe for small sample sizes.

Let $a_0=0$ and $a_k=3\{3L^2+3L+2-6(2L+1)k+10k^2\}/\{L(L-1)(L-2)\}$
for $k=1,\dots,L$.
By Lemma A5 in \cite{Muller1999}, we have $\sum_{k=1}^L a_k =1$.
Then $\hat\sigma^2_{\mathrm{MS}}$ can be represented as the
quadratic form, $\hat\sigma^2_{\mathrm{MS}}=\bfy^T M\bfy$,
where $M=(m_{ij})_{n\times n}$ is a symmetric matrix with elements
\begin{eqnarray*}
m_{ij}=\cases{\displaystyle
1 + \sum_{k=0}^{i-1} a_k, &\quad $i=j=1,\dots,L$, \vspace*{2pt}\cr
2, &\quad $i=j=L+1,\dots,n-L$, \vspace*{2pt}\cr
\displaystyle\sum_{k=i}^n a_{k+L-n}, &\quad $i=j=n-L+1,\dots,n$, \vspace
*{2pt}\cr
-a_{j-i}, &\quad $0<j-i\le L$ and $i\le n-L$, \vspace*{2pt}\cr
-a_{i-j}, &\quad $0<i-j\le L$ and $j\le n-L$, \vspace*{2pt}\cr
0, &\quad otherwise.
}
\end{eqnarray*}
Let $\operatorname{diag}(M)$ denote the diagonal matrix of $M$. By
Dette, Munk and Wagner \cite
{Dette1998} we have
\begin{eqnarray}\label{mse}
\operatorname{MSE}(\hat\sigma^2_{\mathrm{MS}})
&=& \bigl[ (\bfg^T M \bfg)^2 + 4\sigma^2 \bfg^T M^2\bfg+ 4\bfg^T
M \operatorname{diag}(M)\mathbf{1} \sigma^3 \gamma_3 \nonumber\\[-8pt]\\[-8pt]
&&\hspace*{4pt}{}+ \sigma^4 \operatorname{tr}[\{\operatorname{diag}(M)\}
^2] (\gamma_4-3) +
2\sigma^4 \operatorname{tr}(M^2)\bigr] /\operatorname{tr}(M)^2,\nonumber
\end{eqnarray}
where the first term in (\ref{mse}) is the squared bias and the last
four terms make up the variance.

\begin{theorem} \label{theorem3}
Assume that $g$ has a bounded second derivative. Then for the equally
spaced design with $n \to\infty$, $L\to\infty$ and $L/n \to0$, we
have the following bias, variance, and the mean squared error for the
estimator (\ref{est.muller}),
\begin{eqnarray}\label{mse.muller}
\operatorname{Bias}(\hat\sigma^2_{\mathrm{MS}}) &=& \mathrm{o}\biggl(
{L^2 \over n^2}
\biggr), \nonumber\\
\operatorname{var}(\hat\sigma^2_{\mathrm{MS}}) &=& \frac
{1}{n}\operatorname{var}(\varepsilon^2) + \frac{73L}{70n^2}\operatorname
{var}(\varepsilon^2) +
\frac{9}{Ln}\sigma^4 + \mathrm{o}\biggl( \frac{L}{n^2} \biggr) + \mathrm
{o}\biggl(
\frac{1}{
Ln} \biggr), \\
\operatorname{MSE}(\hat\sigma^2_{\mathrm{MS}}) &=& \frac
{1}{n}\operatorname{var}(\varepsilon^2) + \frac{73L}{70n^2}\operatorname
{var}(\varepsilon^2) +
\frac{9}{Ln}\sigma^4 + \mathrm{o}\biggl( \frac{L}{n^2} \biggr) + \mathrm
{o} \biggl(
\frac{1}{
Ln} \biggr) + \mathrm{o}\biggl( \frac{L^4}{n^4} \biggr).\nonumber
\end{eqnarray}
\end{theorem}

Proof of Theorem \ref{theorem3} can be found in Appendix \ref{appD}.
The asymptotical optimal bandwidth is
$L_{\mathrm{opt}}=\sqrt{630n\sigma^4/73\operatorname{var}(\varepsilon
^2)}$. Substituting
$L_{\mathrm{opt}}$ into (\ref{mse.muller}) leads to
\begin{eqnarray}\label{mse.opt.muller}
\operatorname{MSE}(\hat\sigma_{\mathrm{MS}}^2 (L_{\mathrm{opt}})) =
\frac{1}{n}
\operatorname{var}(\epsilon^2) + \frac{\sqrt{45990}}{35}\{\sigma^4
\operatorname{var}(\varepsilon^2)\}^{1/2} n^{-3/2} + \mathrm{o}( n^{-3/2}
).
\end{eqnarray}
The optimal MSE of $\hat\sigma^2$ is \cite{Tong2005}
\[
\operatorname{MSE}(\hat\sigma^2 (m_{\mathrm{opt}}))= \frac{1}{n}
\operatorname{var}(\epsilon^2)
+ \frac{\sqrt{567}}{28} \{ \sigma^4 \operatorname{var}(\varepsilon^2)\}
^{1/2} n^{-3/2} + \mathrm{o}( n^{-3/2}
).
\]
It is clear that both $\hat\sigma^2$ and
$\hat\sigma^2_{\mathrm{MS}}$ reach the optimal efficiency bound
with the same first order term. However, the coefficient of the higher
order term for $\hat\sigma^2_{\mathrm{MS}}$ is about seven times
of that for~$\hat\sigma^2$. Since the higher order term is not
negligible for small to moderate sample sizes, $\hat\sigma^2$ often
provides a much smaller MSE than $\hat\sigma^2_{\mathrm{MS}}$ in
such situations. See simulation results in Section \ref{sec6}.

Even though the two estimators $\hat\sigma^2$ and
$\hat\sigma^2_{\mathrm{MS}}$ look similar for one-dimensional
equally spaced case, there is a fundamental difference behind the
motivations for these estimators: the regression estimator in
\cite{Tong2005} was developed to estimate variances in
nonparametric regression on general domains while the regression
estimator in \cite{Muller1999} was developed for assessing
whether a one-dimensional mean function is smooth. Specifically,
consider model \eqref{model} with $x_i \in{\mathcal{T}}$ where
${\mathcal{T}}$ is an
arbitrary subset in a normed space. Let $d_{ij}=\|x_i-x_j\|^2$
and $s_{ij}={1\over2} (y_i-y_j)^2$ for all pairs $i$ and $j$,
where $1 \leq i <j \leq n$. We fit the following simple linear model
\begin{equation}\label{gest}
s_{ij}= \beta_0 + \beta_1 d_{ij} + \epsilon_{ij},
\qquad d_{ij} \leq m,
\end{equation}
using the least squares where $m>0$ is the bandwidth. The estimate
of $\sigma^2$ is $\hat\sigma^2=\hat\beta_0$. The
variance estimator in \cite{Muller1999} requires an ordering of
the design points which may not be available for a general domain.

%
\begin{table}
\def\arraystratch{0.9}
\caption{Relative mean squared errors for the two estimators with
bandwidths $m_s=L_s=n^{1/2}$ and $m_t=L_t=n^{1/3}$, respectively}\label{tab1}
 \vspace*{-3pt}
\begin{tabular*}{\textwidth}{@{\extracolsep{\fill}}lllllll@{}}
\hline
$n$ & $\sigma^2$ & $g$ & $\hat\sigma^2(m_s)$ & $\hat\sigma^2(m_t)$ &
$\hat\sigma^2_{\mathrm{MS}}(L_s)$ & \multicolumn{1}{c@{}}{$\hat\sigma
^2_{\mathrm{MS}}(L_t)$} \\
\hline
\phantom{00}30 & 0.25 & $g_1$ & 1.33 & 1.58 & 3.97 & 10.80 \\
& & $g_2$ & 1.34 & 1.57 & 3.97 & 10.79 \\
& & $g_3$ & 8.64 & 2.19 & 6.91 & 11.60 \\[2pt]
& 4\phantom{00.} & $g_1$ & 1.32 & 1.57 & 3.91 & 10.75 \\
& & $g_2$ & 1.32 & 1.57 & 3.91 & 10.75 \\
& & $g_3$ & 1.38 & 1.59 & 4.02 & 10.83 \\ [5pt]
\phantom{0}100 & 0.25 & $g_1$ & 1.25 & 1.43 & 2.09 & \phantom{0}5.53 \\
& & $g_2$ & 1.25 & 1.43 & 2.08 & \phantom{0}5.55 \\
& & $g_3$ & 2.06 & 1.45 & 2.30 & \phantom{0}5.50 \\[2pt]
& 4\phantom{00.} & $g_1$ & 1.25 & 1.43 & 2.09 & \phantom{0}5.54 \\
& & $g_2$ & 1.25 & 1.43 & 2.08 & \phantom{0}5.54 \\
& & $g_3$ & 1.27 & 1.43 & 2.09 & \phantom{0}5.52 \\ [5pt]
1000 & 0.25 & $g_1$ & 1.18 & 1.30 & 1.35 & \phantom{0}1.83 \\
& & $g_2$ & 1.18 & 1.30 & 1.35 & \phantom{0}1.83 \\
& & $g_3$ & 1.19 & 1.30 & 1.35 & \phantom{0}1.83 \\[2pt]
& 4\phantom{00.} & $g_1$ & 1.18 & 1.30 & 1.35 & \phantom{0}1.83 \\
& & $g_2$ & 1.18 & 1.30 & 1.35 & \phantom{0}1.83 \\
& & $g_3$ & 1.18 & 1.30 & 1.35 & \phantom{0}1.83 \\
\hline
\end{tabular*}
\end{table}

For the purpose of
illustration, consider the Lake Acidity Data which contains
measurements of 112 lakes in the southern Blue Ridge mountains area
\cite{GuWahba93a}. Of interest is the dependence of the water pH level
($ph$) on the calcium concentration in $\log_{10}$ milligrams per
liter ($t_1$) and the geographical location ($\mathbf{t}_2=(t_{21},t_{22})$
with $t_{21}$ $=$ latitude and $t_{22}$ $=$ longitude). For illustration,
we consider the nonparametric regression model
(\ref{model}) with three different cases of $x$: $x=t_1$, $x=\mathbf{
t}_2$ and $x=(t_1, \mathbf{t}_2)$. These three cases correspond to
three different domains of one, two and three dimensions,
respectively. For the first two cases, we use simple Euclidean
norms. For the third case, we rescale $t_1$ and $\|\mathbf{t}_2\|$ to the
same scale before estimating the variance. Estimates of $\sigma^2$
for the above three cases with $m=n^{1/2}$ are $0.0821$, $0.0884$
and $0.0544$, respectively, using our method. The method in
\cite{Muller1999} does not apply to any one of these three
cases.

\section{Simulation studies} \label{sec6}

In this section, we conduct simulations to compare the performance of
the estimators $\hat\sigma^2$ and $\hat\sigma^2_{\mathrm{MS}}$.
The design points are $x_i=i/n$ and $\varepsilon_i$ are independent
and identically distributed from $N(0,\sigma^2)$.
We consider three mean functions,
$g_1(x)=5x$, $g_2(x)=5x(1-x)$ and $g_3(x)=5\sin(2\uppi x)$.
Note that the first two functions were used in \cite{Muller1999}
and the last one was used in \cite{Tong2005}. We set
coefficients of all three functions to be $5$.
For each mean function, we consider $n=30$, $100$ and $1000$,
corresponding to small, moderate and large sample sizes respectively,
and $\sigma^2=0.25$ and~$4$, corresponding to small and large
variances, respectively.
In total, we have 18 combinations of simulation settings.

For each simulation setting, we generate observations and compute the
estimators $\hat\sigma^2(m)$ and~$\hat\sigma^2_{\mathrm{MS}}(L)$.
For the bandwidth $m$, we choose $m_s=n^{1/2}$ and $m_t=n^{1/3}$ as
suggested in \cite{Tong2005}. For the bandwidth $L$,
M{\"u}ller and Stadtm{\"u}ller \cite{Muller1999} observed that the estimator
$\hat\sigma^2_{\mathrm{MS}}$ is quite stable and does not vary much
with $L$. Therefore, we also choose $L_s=n^{1/2}$ and $L_t=n^{1/3}$
for ease of comparison. The cross-validation method may also be used
to select the bandwidth $m$ in $\hat\sigma^2(m)$ \cite{Tong2005}.
Nevertheless, we did not include this option in our simulations since
the cross-validation method is not readily available for the estimator
$\hat\sigma^2_{\mathrm{MS}}$.

We repeat the simulation 1000 times and
compute the relative mean squared errors\linebreak[4] $n\operatorname{MSE}/(2\sigma^4)$.
Table \ref{tab1} lists relative mean squared errors for all simulation
settings. Note that neither $D$ nor $M$ is guaranteed to
be positive definite. Therefore, $\hat\sigma^2$ and
$\hat\sigma^2_{\mathrm{MS}}$ may take negative values.
Simulations indicate that a negative estimate occurs very rarely for
$\hat\sigma^2$ \cite{Tong2005}, while $\hat\sigma^2_{\mathrm{MS}}$
tends to be negative when $L$ is large \cite{Muller1999}.
We replace negative estimates by zero in the calculation of the
relative mean squared errors.

%
\begin{figure}

\includegraphics{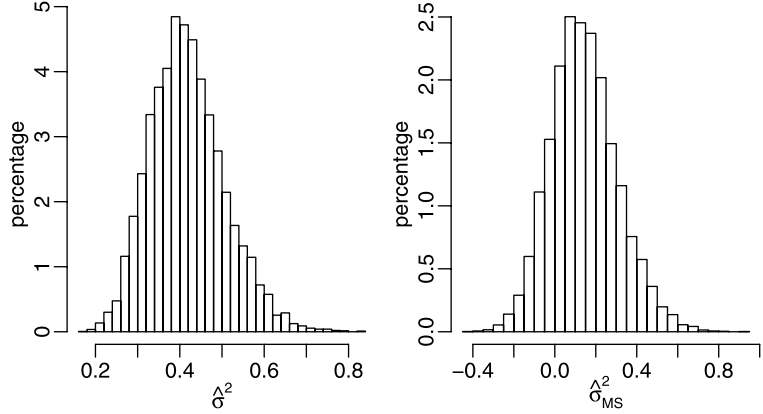}

\caption{Histograms of the variance estimates $\hat\sigma^2(m_s)$
(left) and $\hat\sigma^2_{\mathrm{MS}}(L_s)$ (right) for the case
$(n,\sigma^2,g)=(30,0.25,g_3)$.}
\label{figure1}
\end{figure}

We observe that $\hat\sigma^2$ has smaller relative mean squared errors
than $\hat\sigma^2_{\mathrm{MS}}$ for all settings except for
the case $(n,\sigma^2,g)=(30,0.25,g_3)$.
For this exceptional case, we plot in Figure~\ref{figure1} the
histograms of the nontruncated estimates (including negative
estimates) $\hat\sigma^2(m_s)$ and $\hat\sigma^2_{\mathrm{MS}}(L_s)$.
A relatively large portion of $\hat\sigma^2_{\mathrm{MS}}(L_s)$
takes negative values.
The choice of the bandwidth $m_s$ is too large for $\hat\sigma^2$
when $n$ is small
\cite{Tong2005}. Overall, the estimator $\hat\sigma^2$ performs
better than $\hat\sigma^2_{\mathrm{MS}}$, confirming the theoretical
results in Section \ref{sec5}.
Comparisons between $\hat\sigma^2(m_s)$ and $\hat\sigma^2(m_t)$ are
similar to those in~\cite{Tong2005}.

\begin{appendix}
\section{\texorpdfstring{Proof of Theorem \protect\ref{theorem1}}{Proof of Theorem 1}}\label{appA}

For ease of notation, let $g_i=g(x_i)$, $i=1,\dots,n$. Write $s_k$ as a
sum of three parts, $s_k=L_1 + L_2 + L_3$, where
\begin{eqnarray*}
L_1 &=& \frac{1}{2(n-k)}\sum_{i=k+1}^n (g_i-g_{i-k})^2, \\
L_2 &=& \frac{1}{n-k}\sum_{i=k+1}^n (g_i-g_{i-k})(\varepsilon
_i-\varepsilon_{i-k}), \\
L_3 &=& \frac{1}{2(n-k)}\sum_{i=k+1}^n (\varepsilon_i-\varepsilon_{i-k})^2.
\end{eqnarray*}
Applying the Taylor expansion, it can be
shown that $L_1 = (k^2/n^2)J + \mathrm{o}({k^2/n^2}) = \mathrm
{o}_p(n^{-1/2})$ when $k=n^r$ with $0<r<3/4$.
For $L_2$, we have
\begin{eqnarray*}
\mathrm{E}(L_2^2)
= \frac{2\sigma^2}{(n-k)^2} \Biggl\{\sum_{i=k+1}^n (g_i-g_{i-k})^2 -
\sum_{i=k+1}^{n-k} (g_i-g_{i-k})(g_{i+1}-g_i) \Biggr\}
= \mathrm{O} \biggl(\frac{k^2}{n^3} \biggr).
\end{eqnarray*}
This implies that $L_2= \mathrm{o}_p(n^{-1/2})$ for any $k=\mathrm{o}(n)$\vspace*{1.5pt}.

Rewrite $L_3$ as $L_3 = \sigma^2 + \sum_{i=k+1}^n \xi_i(k)/(n-k)$,
where $\xi_i(k) = (\varepsilon_i-\varepsilon_{i-k})^2/2 -\sigma^2$.
For any given $k$, $\{\xi_i(k), i=k+1,\dots,n\}$ is a strictly
stationary sequence of random variables
with mean zero and autocovariance function
\begin{eqnarray*}
\gamma(\tau)
= \gamma(s,s+\tau)
= \cases{
(\gamma_4+1)\sigma^4/2, &\quad $ \tau=0$, \vspace*{2pt}\cr
(\gamma_4-1)\sigma^4/4, &\quad $\tau=k$, \vspace*{2pt}\cr
0, &\quad otherwise.
}
\end{eqnarray*}
Note also that the sequence $\{\xi_i(k), i=k+1,\dots,n\}$ is
$m$-dependent with $m=k$.
Thus by the central limit theorem for strictly stationary $m$-dependent
sequences \cite{Brockwell1991},
$\sqrt{n}(L_3 -\sigma^2) \stackrel{\mathcal{D}}\to N(0,\nu_k^2)$ as
$n\to\infty$,
where $\nu_k^2 = \gamma(0) + 2\sum_{\tau=1}^k \gamma(\tau) = \gamma
_4\sigma^4$.
Finally, noting that $s_k = L_1+L_2+L_3 = L_3 + \mathrm{o}_p(n^{-1/2})$,
we have $\sqrt{n}(s_k -\sigma^2) \stackrel{\mathcal{D}}\to N(0, \gamma
_4\sigma^4)$ as $n\to\infty$.

\section{\texorpdfstring{Proof of Theorem \protect\ref{theorem2}}{Proof of Theorem 2}}\label{appB}

We first state two lemmas. Lemma \ref{lemma1} is an immediate result
from \cite{Whittle1964}. Lemma \ref{lemma2} was derived, in
essence, in \cite{Tong2005}.

\begin{lemma} \label{lemma1}
Assume that the matrix $A=(a_{ij})_{n\times n}$ satisfies
$a_{ij}=a_{i-j}$ and $\sum_{-\infty}^\infty a_k^2 <
\infty$. Furthermore, assume that $\mathrm{E}(\varepsilon^6)$ is finite.
Then
\begin{eqnarray*}
\frac{1}{n}\bfvarepsilon^TA \bfvarepsilon= \frac{1}{n}\sum_{i=1}^n\sum
_{j=1}^n a_{i-j}\varepsilon_i\varepsilon_j
\stackrel{\mathcal{D}}\longrightarrow N(a_0\sigma^2, \sigma
_A^2),\qquad \mbox{as } n\to\infty,
\end{eqnarray*}
where $\sigma_A^2 = (\gamma_4-3)a_0^2\sigma^4/n + 2\sigma^4\sum
_{i=1}^n\sum_{j=1}^n a_{i-j}^2/n^2$.
\end{lemma}

\begin{lemma} \label{lemma2}
Assume that $m \to\infty$ and $m/n \to0$. Then
\begin{longlist}[(iii)]
\item[(i)] $\sum_{k=1}^m b_k = m - \frac{5m^2}{16n}+ \mathrm{o}(m)$;
\item[(ii)] $\sum_{k=j}^m b_k = m - \frac{9}{4}j + \frac{5j^3}{4m^2}
+ \mathrm{o}(m), 1 \leq j \leq m$;
\item[(iii)] $\sum_{k=1}^m b_k^2 = \frac{9}{4}m + \mathrm{o}(m)$;
\item[(iv)] $\bfg^T D\bfg= \mathrm{O}(m^4/n^2)$;
\item[(v)] $\bfg^T D^2 \bfg= \mathrm{O}(m^5/n^2)$.
\end{longlist}
\end{lemma}

\renewcommand{\theequation}{\arabic{equation}}
\setcounter{equation}{12}
\begin{pf*}{Proof of Theorem \ref{theorem2}}
Noting that $\bfy=\bfg+\bfvarepsilon$ and $\operatorname{tr}(D)=2N$, we have
\begin{eqnarray}\label{approx1}
\hat\sigma^2 = \frac{1}{2N}\bfg^T D\bfg+ \frac{1}{N}\bfg^T
D\bfvarepsilon+ \frac{1}{2N}\bfvarepsilon^T D\bfvarepsilon.
\end{eqnarray}
The first term in (\ref{approx1}) corresponds to the bias term of the
least squares estimator.
By Lemma \ref{lemma2}, we have $\bfg^T D\bfg/(2N) = \mathrm{O}(m^3/n^3)$.
Thus, for any $m=n^r$ with $0<r<5/6$,
\begin{eqnarray}\label{term1}
\frac{1}{2N} \bfg^T D\bfg= \mathrm{o}(n^{-1/2}).
\end{eqnarray}
For the second term in (\ref{approx1}), by Lemma \ref{lemma2} we have
$\mathrm{E} (\bfg^T D\bfvarepsilon/N )^2 = \bfg^T D^2 \bfg/N^2 = \mathrm{O}({m^3/n^4})$.
This implies that, for any $m=\mathrm{o}(n)$,
\begin{eqnarray}\label{term2}
\frac{1}{N}\bfg^T D\bfvarepsilon= \mathrm{o}_p(n^{-1/2}).
\end{eqnarray}

Now we derive the limiting distribution of the third term in (\ref{approx1}).
Let $nD/(2N)=C - H$, where $C = (c_{ij})_{n\times n}$ with elements
\begin{eqnarray*}
c_{ij}= \cases{\displaystyle
n\sum_{k=1}^m b_k/N, &\quad $1 \le i=j \le n$, \vspace*{2pt}\cr
-nb_{|i-j|}/(2N), &\quad $0<|i-j| \leq m$, \vspace*{2pt}\cr
0, &\quad  otherwise,
}
\end{eqnarray*}
and $H=\operatorname{diag}(h_1, h_2, \dots, h_n)$ with elements $h_{i}=n\sum
_{\min(i,n+1-i,m+1)}^{m+1} b_k/(2N)$.
Then
\begin{eqnarray}\label{partition1}
\frac{1}{2N}\bfvarepsilon^T D\bfvarepsilon= \frac{1}{n}\bfvarepsilon^T
C\bfvarepsilon- \frac{1}{n}\bfvarepsilon^T H\bfvarepsilon.
\end{eqnarray}
For the matrix $C$, let $c_{ij}=c_{i-j}$ with $c_0=n\sum_{k=1}^m b_k/N$,
$c_{i-j}=c_{j-i}=-nb_{|i-j|}/(2N)$ for $0<|i-j| \le m$, and
$c_{i-j}=c_{j-i}=0$ for $|i-j|>m$.
By Lemma \ref{lemma2}, for any $m=\mathrm{o}(n)$, $\sum_{-\infty
}^\infty c_k^2 = c_0^2 + 2\sum_{k=1}^m c_k^2 = 1 + \mathrm{o}(1)<\infty$.
Then under the assumption that $\mathrm{E}(\varepsilon^6)$ is finite,
by Lemma \ref{lemma1} we have
\begin{eqnarray}\label{term3}
\sqrt{n} \biggl(\frac{1}{n}\bfvarepsilon^T C\bfvarepsilon- c_0\sigma^2 \biggr)
\stackrel{\mathcal{D}}\longrightarrow N(0, \sigma_c^2),\qquad \mbox{as } n\to
\infty,
\end{eqnarray}
where
\[
\sigma_c^2 = \frac{n^2(\gamma_4-1)\sigma^4}{N^2} \Biggl(\sum_{k=1}^m b_k \Biggr)^2 +
\frac{n\sigma^4}{N^2}\sum_{k=1}^m (n-k)b_k^2.
\]
For the second term in (\ref{partition1}), note that $\bfvarepsilon^T
H\bfvarepsilon=\sum_{1}^m h_{i}\varepsilon_i^2 + \sum_{n-m+1}^n
h_{i}\varepsilon_i^2$.
By Lemma \ref{lemma2}, it is easy to see that
\begin{eqnarray*}
\mathrm{E} \Biggl(\sum_{i=1}^m h_{i}\varepsilon_i^2 \Biggr)^2
&=& (\gamma_4-1)\sigma^4 \frac{n^2}{4N^2}\sum_{i=1}^m \Biggl(\sum_{\min
(i,n+1-i,m+1)}^{m+1} b_k \Biggr)^2 \\
&&{} + \frac{n^2\sigma^4}{4N^2} \Biggl(\sum_{i=1}^m \sum_{\min
(i,n+1-i,m+1)}^{m+1} b_k \Biggr)^2 \\
&=& \mathrm{O}(m^2).
\end{eqnarray*}
Similarly, we have $\mathrm{E} (\sum_{n-m+1}^n h_{i}\varepsilon_i^2 )^2=\mathrm{O}(m^2)$.
This leads to $\mathrm{E}(\bfvarepsilon^T H\bfvarepsilon/n)^2=\mathrm{O}(m^2/n^2)$.
Further, for any $m=n^r$ with $0<r<1/2$,
\begin{eqnarray}\label{term4}
\frac{1}{n}\bfvarepsilon^T H\bfvarepsilon=\mathrm{o}_p(n^{-1/2}).
\end{eqnarray}

Combining (\ref{term1}), (\ref{term2}), (\ref{term3}) and
(\ref{term4}), and applying the Slutsky theorem, we have
\begin{eqnarray}\label{appen2.f1}
\frac{\sqrt{n}(\hat\sigma^2 - c_0\sigma^2)}{\sigma_c} \stackrel{\mathcal
{D}}\longrightarrow N(0, 1),\qquad \mbox{as } n\to\infty.
\end{eqnarray}
Note also that, by Lemma \ref{lemma2},
\begin{eqnarray*}
c_0 &=& \frac{n}{nm-m(m+1)/2} \biggl\{m - \frac{5m^2}{16n}+ \mathrm{o}(m) \biggr\} = 1
+ \mathrm{O} \biggl(\frac{m}{n} \biggr), \\
\sigma_c^2 &=& \frac{n^2(\gamma_4-1)\sigma^4}{N^2} \Biggl(\sum_{k=1}^m b_k \Biggr)^2
+ \frac{n\sigma^4}{N^2}\sum_{k=1}^m (n-k)b_k^2 = (\gamma_4-1)\sigma^4 +
\mathrm{o}(1).
\end{eqnarray*}
Thus for any $m=n^r$ with $0<r<1/2$, we have $\sqrt{n}(c_0-1)=\mathrm{o}(1)$.
In addition, $(\gamma_4-1)\sigma^4/\sigma_c^2 \to1$ as $n\to\infty$.
Then by (\ref{appen2.f1}) and the Slutsky theorem,
\begin{eqnarray*}
\frac{\sqrt{n}(\hat\sigma^2 - \sigma^2)}{\sqrt{(\gamma_4-1)\sigma^4}}
&=& \frac{\sigma_c}{\sqrt{(\gamma_4-1)\sigma^4}} \biggl\{\frac{\sqrt{n}(\hat\sigma
^2 - c_0\sigma^2)}{\sigma_c} + \frac{\sqrt{n}(c_0-1)\sigma^2}{\sigma
_c} \biggr\} \\
&\stackrel{\mathcal{D}}\longrightarrow& N(0, 1),\qquad \mbox{as } n\to\infty.
\end{eqnarray*}
\upqed\end{pf*}

\section{Derivation of covariances between Rice estimators}\label{appC}

For any $1\leq b<k=\mathrm{o}(n)$, we have
\begin{eqnarray*}
\mathrm{E}(s_bs_k)
&=& \frac{1}{4(n-b)(n-k)}\\
&&{}\times \Biggl\{\sum_{i=k+1}^n \mathrm
{E}(y_i-y_{i-k})^2(y_{i-k+b}-y_{i-k})^2 + \sum_{i=k+1}^n \mathrm{E}(y_i-y_{i-k})^2(y_i-y_{i-b})^2 \\
&&\hspace*{16pt}{} + \sum_{i=k+b+1}^n \mathrm{E}(y_i-y_{i-k})^2(y_{i-k}-y_{i-k-b})^2 + \sum_{i=k+1}^{n-b} \mathrm{E}(y_i-y_{i-k})^2(y_{i+b}-y_{i})^2 \\
&&\hspace*{16pt} {}+ \sum_{(i,j)\in\mathcal{E}} \mathrm{E}(y_i-y_{i-k})^2(y_j-y_{j-b})^2
\Biggr\} \\
&=& \frac{1}{4(n-b)(n-k)}(I_1+I_2+I_3+I_4+I_5),
\end{eqnarray*}
where $\mathcal{E}=\{(i,j)\dvtx i=k+1,\dots,n; j=b+1,\dots,n; i\neq j; i\neq
j-b; i-k\neq j; i-k\neq j-b\}$.
It is easy to verify that
$I_1 + I_2 = 2(n-k)(\gamma_4+3)\sigma^4 + \mathrm{O}({k^2/n})$,
$I_3 + I_4 = 2(n-k-b)(\gamma_4+3)\sigma^4 + \mathrm{O}({k^2/n})$
and $I_5 = 4 \{(n-k)(n-b)-2(2n-2k-b) \}\sigma^4 + 4\sigma
^2(n-b)(n-k)(b^2+k^2)J/n^2 + \mathrm{O}({k^3/n})$.
Therefore,
\[
\mathrm{E}(s_bs_k) = \frac{2n-2k-b}{2(n-b)(n-k)}(\gamma_4-1)\sigma^4 +
\sigma^4 + \frac{b^2+k^2}{ n^2}J\sigma^2 + \mathrm{O} \biggl(\frac{k^3}{n^3} \biggr).
\]
Note also that $\mathrm{E}(s_b)=\sigma^2 + Jd_b + \mathrm{O}(b^3/n^3) + \mathrm
{o}(1/n^2)$ and $\mathrm{E}(s_k)=\sigma^2 + Jd_k + \mathrm{O}(k^3/n^3) + \mathrm
{o}(1/n^2)$.
Thus,
\begin{eqnarray*}
\operatorname{Cov}(s_b,s_k) = \frac{2n-2k-b}{2(n-b)(n-k)}(\gamma_4-1)\sigma^4 + \mathrm{O}
\biggl(\frac{k^3}{n^3} \biggr) + \mathrm{o} \biggl(\frac{1}{n^2} \biggr).
\end{eqnarray*}
Finally, for any $k=n^r$ with $0<r<2/3$, we have $k^3/n^3=\mathrm{o}(1/n)$
and therefore $\operatorname{Cov}(s_b,s_k) = (\gamma_4-1)\sigma^4/n + \mathrm{o}(1/n)$.

\section{\texorpdfstring{Proof of Theorem \protect\ref{theorem3}}{Proof of Theorem 3}}\label{appD}

\begin{lemma} \label{Muller}
Assume that $g$ has a bounded second derivative. Then for the equally
spaced design with $n \to\infty$, $L\to\infty$ and $L/n\to0$, we have
\begin{enumerate}[(iii)]
\item[(i)] $\operatorname{tr}(M) = 2(n-L)$;
\item[(ii)] $\operatorname{tr}[\{\operatorname{diag}(M)\}^2] = 4n - 134L/35 +
\mathrm{o}(L)$;
\item[(iii)] $\operatorname{tr}(M^2) = 4n - 134L/35 + 18n/L + \mathrm
{o}(L) + \mathrm{o} ( {n/L} )$;
\item[(iv)] $\bfg^T M^2\bfg= \mathrm{O} ( {L^3/n^2} )$;
\item[(v)] $\bfg^T M \operatorname{diag}(M)\mathbf{1} = \mathrm{O} ( {L^2/n} )$.
\end{enumerate}
\end{lemma}

\begin{pf}
It is easy to verify that $\sum_{k=1}^L a_k = 1$, $\sum_{k=1}^i a_k =
{9i/L} - {18i^2/L^2} + {10i^3/L^3} + \mathrm{o}({i/L})$ for $1\leq
i\leq L$,
$\sum_{k=1}^L a_k^2 = {9/L} + \mathrm{o}({1/L})$, $\sum_{k=1}^L ka_k = \mathrm{O}(L)$
and $\sum_{k=1}^L k^2a_k = \mathrm{O}(L^2)$.

(i)
$\operatorname{tr}(M) = 2L\sum_{k=1}^L a_k + 2(n-2L)\sum_{k=1}^L a_k = 2(n-L)$.

(ii)
Note that $a_0=0$ and $\sum_{k=n-L+i}^n a_{k+L-n}= 1 - \sum_{k=0}^{i-1}
a_k$. We have
\begin{eqnarray*}
\operatorname{tr}[\{\operatorname{diag}(M)\}^2]
&=& 4(n-2L) + \sum_{i=1}^L \Biggl(1 + \sum_{k=0}^{i-1} a_k \Biggr)^2 +
\sum_{i=1}^L \Biggl(1 - \sum_{k=0}^{i-1} a_k \Biggr)^2 \\
&=& 4n-6L + 2\sum_{i=1}^L \biggl\{\frac{9i}{ L} - \frac{18i^2}{ L^2} +
\frac{10i^3}{ L^3} + \mathrm{o} \biggl( \frac{i}{ L} \biggr) \biggr\}^2 \\
&=& 4n - \frac{134}{35}L + \mathrm{o}(L).
\end{eqnarray*}

(iii)
By (ii), we have
\begin{eqnarray*}
\operatorname{tr}(M^2)
&=& \operatorname{tr}[\{\operatorname{diag}(M)\}^2] + \sum_{i=1}^L \Biggl(\sum_{k=1}^L a_k^2 +
\sum_{k=0}^{i-1}a_k^2 \Biggr) + 2\sum_{i=L+1}^{n-L}\sum_{k=1}^L a_k^2 + \sum
_{i=1}^L \sum_{k=i}^L a_k^2 \\
&=& \operatorname{tr}[\{\operatorname{diag}(M)\}^2] + 2(n-L)\sum_{k=1}^L a_k^2 \\
&=& 4n - \frac{134}{35}L + \frac{18n}{ L} + \mathrm{o}(L) + \mathrm{o}
\biggl(\frac{n}{ L} \biggr).
\end{eqnarray*}

(iv)
Noting that $M$ is a symmetric matrix, we have $\bfg^T M^2\bfg=
(M\bfg)^T M\bfg\triangleq\bfh^T\bfh$ where
$\bfh=M\bfg=(h_1,\dots,h_n)^T$.
Under the condition that $g$ has a bounded second derivative, it is
easy to verify that for $i\in[L+1,n-L]$,
\begin{eqnarray*}
h_i = \sum_{k=1}^L a_k(g_i-g_{i-k}) - \sum_{k=1}^L a_k (g_{i+k}-g_i)
= - \frac{1}{n^2}g^{\prime\prime}_i \sum_{k=1}^L k^2a_k + \mathrm{o}
\biggl(\frac{m^3}{
n^2} \biggr)
= \mathrm{O} \biggl(\frac{L^2 }{n^2} \biggr).
\end{eqnarray*}
Similarly, we can show that for $i\in[1,L]$ or $i\in[n-L+1,n]$,
$h_i=\mathrm{O}(L/n)$. Finally,
\begin{eqnarray*}
\bfg^T M^2\bfg= \bfh^T \bfh= \sum_{i=1}^L h_i^2 +\sum_{i=L+1}^{n-L}
h_i^2 +\sum_{i=n-L+1}^n h_i^2 =\mathrm{O} \biggl(\frac{L^3}{ n^2} \biggr).
\end{eqnarray*}

(v)
Note that $\bfg^T [M \operatorname{diag}(M)\mathbf{1}] = (M\bfg)^T \operatorname{diag}(M)\mathbf{1} = \bfh^T \operatorname{diag}(M)\mathbf{1}$. We have
\begin{eqnarray*}
\bfg^T [M \operatorname{diag}(M)\mathbf{1}]
= \sum_{i=1}^L h_i \cdot \mathrm{O}(1) + \sum_{i=L+1}^{n-L} h_i \cdot \mathrm{O}(1) +
\sum_{i=n-L+1}^n h_i \cdot \mathrm{O}(1)
= \mathrm{O} \biggl(\frac{L^2 }{ n} \biggr).
\end{eqnarray*}
\upqed\end{pf}

\begin{pf*}{Proof of Theorem \ref{theorem3}}
By M{\"u}ller and Stadtm{\"u}ller \cite{Muller1999}, $\operatorname{Bias}(\hat\sigma^2_{\mathrm{MS}})=\bfg
^T M \bfg/\operatorname{tr}(M)=\mathrm{o}(L^2/n^2)$.
Note that the last four terms in (\ref{mse}) make up the variance.
By Lemma \ref{Muller} and the facts that $L/n\to0$ and $\sigma^4
(\gamma_4-3)= \operatorname{var}(\varepsilon^2) - 2\sigma^4$, we have
\begin{eqnarray*}
\operatorname{var}(\hat\sigma^2_{\mathrm{MS}})
&=& \frac{1}{4(n-L)^2} \biggl[ \{\operatorname{var}(\varepsilon^2)-2\sigma^4 \} \biggl\{4n -
\frac{134}{35}L + \mathrm{o}(L) \biggr\} \\
&&{} + 2\sigma^4 \biggl\{4n - \frac{134}{35}L + \frac{18n}{ L} + \mathrm{o}(L) +
\mathrm{o} \biggl(\frac{n}{ L} \biggr) \biggr\} \biggr] \\
&=& \frac{1}{4(n-L)^2} \biggl\{ \biggl(4n - \frac{134}{35}L \biggr)\operatorname{var}(\varepsilon
^2) + \frac{36n}{L}\sigma^4 + \mathrm{o}(L) + \mathrm{o} \biggl(\frac{n}{ L} \biggr)
\biggr\} \\
&=& \frac{1}{ n}\operatorname{var}(\varepsilon^2) + \frac{73L}{70n^2}\operatorname{var}(\varepsilon^2) + \frac{9}{ Ln}\sigma^4 +
\mathrm{o} \biggl(\frac{L}{ n^2} \biggr)
+ \mathrm{o} \biggl(\frac{1}{ Ln} \biggr).
\end{eqnarray*}
Finally, we have
\begin{eqnarray*}
\operatorname{MSE}(\hat\sigma^2_{\mathrm{MS}})
&=& \frac{1}{ n}\operatorname{var}(\varepsilon^2) + \frac{73L}{70n^2}\operatorname{var}(\varepsilon^2) + \frac{9}{Ln}\sigma^4 + \mathrm{o} \biggl(\frac{L}{n^2} \biggr)
+ \mathrm{o} \biggl(\frac{1}{ Ln} \biggr) + \mathrm{o} \biggl(\frac{L^4}{ n^4} \biggr).
\end{eqnarray*}
\upqed\end{pf*}
\end{appendix}

\section*{Acknowledgements}

Tiejun Tong's research was supported by Hong Kong RGC Grant HKBU202711,
and Hong Kong Baptist University Grants FRG1/10-11/031 and FRG2/10-11/020.
Yanyuan Ma's research was supported by NSF Grant DMS-09-06341 and NINDS
Grant R01-NS073671.
Yuedong Wang's research was supported by NSF Grant DMS-07-06886.
The authors thank the editor, the associate editor, and a referee for
their constructive comments that substantially improved an earlier draft.

%

\printhistory

\end{document}